# Density of the Values Set of the Tau Function
N. A. Carella

**Abstract**: It is shown that the density of the values set $\{\tau(n) : n \leq x\}$ of the $n$th coefficients $\tau(n) \in \mathbb{Z}$ of the discriminant function $\Delta(z) = \sum_{n=1}^{\infty} \tau(n) q^n$, a cusp form of level $N = 1$ and weight $k = 12$, has the lower bound $\#\{\tau(n) : n \leq x\} \gg x/\log x$. The known density is $\#\{\tau(n) : n \leq x\} \gg x^{1/2+o(1)}$, and the expected density is $\#\{\tau(n) : n \leq x\} \gg x$. The solutions set of the equation $\tau(p) = 0$ for all primes $p \geq 2$, which arises as a singular case of this analysis, is discussed within.

## 1. Introduction
Let $\tau_k(n)$ be the $n$th coefficient of the discriminant function $\Delta_k(z) = \sum_{n=1}^{\infty} \tau_k(n) q^n$, a cusp form of level $N = 1$ and weight $k \geq 12$, confer [1, Chapter 3], [9, p. 20], [22, p. 97], [30, p. 13] for other details. The mutiltiplicative property $\tau_k(mn) = \tau_k(m) \tau_k(n)$, $\gcd(m, n) = 1$, see [11], [22, p. 105], indicates that subset of prime integers $\{p : \tau_k(p) \neq 0\}$ is a multiplicative basis of the image $\tau_k(\mathbb{N}) = \{\tau_k(n): n \in \mathbb{Z}\} \subset \mathbb{Z}$ of the $n$th coefficients $\tau_k(n)$. The multiplicative span of the primes basis is believed to satisfy the asymptotic formula $\#\{\tau_k(n) : n \leq x\} \sim x$. It is also believed that the correspondence $\tau_k(n) \longrightarrow \mathbb{Z}$ is injective.

For the first case $\tau(n) = \tau_{12}(n)$, the lower estimate $\#\{\tau(n) : n \leq x\} \gg x^{1/3+o(1)}$ of the density of the value set of the $n$th coefficients $\tau(n) \in \mathbb{Z}$ of the discriminant function $\Delta(z) = \sum_{n=1}^{\infty} \tau(n) q^n$, was proved in [28]. And in [7] it is shown that it has a lower bound

$$\#\{\tau(n) : n \leq x\} \gg x^{1/2} e^{-4 \log x/\log\log x}. \tag{1}$$

Conditional results of the form $\#\{\tau_k(n) : n \leq x\} \geq .999\, x$, $k = 12, 16, 18, 20, 22$, and $26$, are proved in [24]. This note improves the lower unconditional estimate (1) of the value set to the following.

**Theorem 1.** The density of the values of the $n$th coefficients $\tau(n)$ of the discriminant function $\Delta(z) = \sum_{n=1}^{\infty} \tau(n) q^n$ satisfies the lower bound

$$\#\{\tau(n) : n \leq x\} \gg x/\log x. \tag{2}$$

The proof of Theorem 1 appears in Section 3. This result is an immediate consequence of Theorems 2 on the finite cardinalities of the subsets of integers $V_t = \{p \in \mathbb{P} : \tau(p) = t\}$, where $t \in \mathbb{Z} = \{0, \pm 1, \pm 2, \pm 3, \ldots\}$ is a fixed integer, and $\mathbb{P} = \{2, 3, 5, \ldots\}$ is the set of prime numbers, which is proved in Section 3. The Lang-Trotter conjecture predicts that the cardinality of any of these subsets associated with cusp forms of level $N = 1$ and weights $k \geq 4$, are finite, see [10], [16], [23]. This is also implied by the Atkin-Serre conjecture $|\tau(p)| \gg p^{(k-3)/2-\epsilon}$ with $\epsilon > 0$ arbitrarily small, see [17], and [23]. Section 5 covers the singular case $\tau(p) = t = 0$. Section 6 covers the solutions set of the equation $\tau(p) = p\, t,\ t \in \mathbb{Z}$.

For the other modular forms of levels $N \geq 1$ and weights $k \leq 3$, the subsets $V_t$ can have infinite cardinalities, for example, the work in [4] for modular forms of weight $k = 2$ and the corresponding subsets of supersingular primes. Other cases are discussed in [24].

The nonvanishing results for the Taylor coefficients $\tau_{z_0}(n) \neq 0$ of the expansions of the discriminant cusp form $\Delta_{z_0}(z) = \sum_{n=1}^{\infty} \tau_{z_0}(n) q^n$ at certain complex points $z_0 \in \mathbb{H} = \{ z \in \mathbb{C} : \mathfrak{Re}(z) > 0 \}$ of the upper half plane are proved in [18]. In Section 4, it is shown that equation $\tau(p) = 0$ for all primes $p \geq 2$ arises as a singular case of the analysis for the case $\tau(p) = t \neq 0$. The usual discriminant function $\Delta(z) = \sum_{n=1}^{\infty} \tau(n) q^n$ is an expansion at the point at infinity $z_0 = \infty$.

## 2. Congruences

Given a fixed integer $m \in \mathbb{Z}$, the congruences of the tau function can be used to derive Diophantine equations and criteria for the existence of integers solutions $n \in \mathbb{N} = \{ 0, 1, 2, 3, ... \}$ of the equation

$$\tau(n) = m. \tag{3}$$

Many of the complex congruences of the tau function are expressed in terms of convolutions of the sum of divisors function $\sigma_s(n) = \sum_{d \mid n} d^s$, $s \geq 0$, For an integer $n \geq 1$, and the parameters $a, b, c, d \in \mathbb{N} = \{ 0, 1, 2, 3, ... \}$, the weighted convolution of the sum of divisors function is defined by

$$S_n(a, b, c, d) = \sum_{k=1}^{n-1} k^a (n-k)^b \sigma_c(k) \sigma_d(n-k). \tag{4}$$

A few of the congruences from the wide range of congruences for the tau function are the Niebur formula

$$\tau(n) = n^4 \sigma(n) - 24 \left( 35 S_n(4, 0, 1, 1) - 52 n S_n(3, 0, 1, 1) + 18 n^2 S_n(2, 0, 1, 1) \right), \tag{5}$$

the Lanphier formula

$$2 n \tau(n) = -n^3 \sigma_7(n) + 3 n^3 \sigma_3(n) + 360 S_n(3, 0, 3, 3), \tag{6}$$

and the more recently developed formulas such as the convolution identities

$$\begin{aligned}
2160 S_n (2, 0, 3, 3) &= 5 \sigma_7(n) - 9 n^2 \sigma_3(n) + 4 \tau(n), \\
1080 S_n (2, 1, 3, 3) &= n^3 \sigma_7(n) - n \tau(n), \\
30888 S_n (2, 2, 3, 3) &= 13 n^4 \sigma_7(n) - 22 n^2 \tau(n) + 9 \tau(n) + 2160 r(n),
\end{aligned} \tag{7}$$

where the last term is defined by the finite sum $r(n) = \sum_{k=1}^{n-1} \sigma_3(k) \tau(n-k)$, see [19, p. 103]. The proofs of many other congruences are given in [3], [20], [21], [22], [5], [6], [31], [32], and similar references.

The simpler congruences such as $\tau(n) \equiv n \sigma(n) \mod 30$, see [3, p. 24], do not offer sufficient flexibility to develop interesting criteria for equation (3), but certain combinations of other complex congruences as (5) to (6) do produce complex criteria suitable to investigate the solutions set of equation (3).

### 3. The Subsets $V_t$

Let $m \in \mathbb{Z}$ be a fixed integer, and let $\beta : \mathbb{N} \longrightarrow \mathbb{Z}$ be an arithmetic function. The valence or multiplicity of the arithmetic function $\beta$ at $m$ is defined by the cardinality of the subset of integers

$$V_m = \{ n \in \mathbb{N} : \beta(n) = m \}. \tag{8}$$

One of the earliest problem on the valence of a function is the Carmichel conjecture for the Euler totient function $\varphi : \mathbb{N} \longrightarrow \mathbb{Z}$. This conjecture claims that $\#V_m = 0$ or $\#V_m = \#\{ n \in \mathbb{N} : \varphi(n) = m \} > 1$ for all fixed integer $m \geq 2$. There are partial results on this problem available in the literature. In this case, the associated Dirichlet series $\sum_{n \geq 1} \varphi(n) n^{-s} = \zeta(s-1)/\zeta(2s)$ is not derived from a cusp form. Exempli gratia, $f(z) = \sum_{n=1}^{\infty} \varphi(n) q^n$.

In the specific case of the $n$th coefficient function $\tau(n)$ of the discriminant function $\Delta(z)$, about two decades ago it proved that the subset of values $V_{2t+1} = \{ n \in \mathbb{N} : \tau(n) = 2t + 1 \}$ is a finite subset, see [17]. And very recently, a conditional result proved that the subset of values $V_t = \{ n \in \mathbb{N} : \tau(n) = t \}$ is a finite subset for any fixed integer $0 \neq t \in \mathbb{Z}$, confer [8]. The $n$th coefficient function $\tau(n) \longrightarrow \mathbb{Z}$ is expected to be one-to-one for all nonzero integers $t \neq 0$, but a subset of zero density. For example, $\#V_t = \{ n \in \mathbb{N} : \tau(n) = t \} = 1$ for almost all nonzero integers $t \neq 0$. But that could be difficult to prove. Here it will be demonstrated that for any fixed ineger $0 \neq t \in \mathbb{Z}$, the subset of values $V_t = \{ n \in \mathbb{N} : \tau(n) = t \}$ is a finite subset unconditionally.

In light of the the multiplicative property $\tau(p q) = \tau(p) \tau(q), \gcd(p, q) = 1$, see [11], [22, p. 105], it is immediate that the first occurrence of any vanishing coefficient $\tau(n) = 0$ occurs at a prime argument $n = p \geq 2$, see [12, Theorem 2].

**Theorem 2.** Let $\tau(p)$ be the $p$th coefficient of the discriminat function $\Delta(z) = \sum_{n=1}^{\infty} \tau(n) q^n$. Then, for each fixed $0 \neq t \in \mathbb{Z}$, the cardinality of the subset $V_t$ is finite: $\#V_t = \#\{ p \in \mathbb{N} : \tau(p) = t \} = O(1)$.

**Proof:** Fix an integer $0 \neq t \in \mathbb{Z}$, and let $p \geq 2$ be a prime variable. To show that the equation $\tau(p) = t \neq 0$ has a finite number of prime solutions, consider the system of congruences

$$\begin{aligned} n^3 \sigma_7(n) - n \tau(n) &\equiv 0 \mod 2^3 \cdot 3^3 \cdot 5, \\ 13 n^4 \sigma_7(n) - 22 n^2 \tau(n) + 9 \tau(n) &\equiv 0 \mod 2^3 \cdot 3^3, \end{aligned} \tag{9}$$

refer to the identities (7) and [19, p. 103] for a proof. Expanding the divisors function $\sigma_s(n) = \sum_{d \mid n} d^s$ at the prime variable $p \geq 2$ produces the system of Diophantine equations

$$\begin{aligned} p^3\left(p^7+1\right)-p\,\tau(p)-4\,X &= 0, \\ 13\,p^4\left(p^7+1\right)-\left(22\,p^2-9\right)\tau(p)-4\,Y &= 0, \\ \tau(p) &= t, \end{aligned} \qquad (10)$$

where $X, Y \in \mathbb{Z}$ are integers, and $p \in \mathbb{P}$ is a prime variable. Any combination of the moduli specified by the system of congruences (9) can be selected; the choice of 4 for both congruences in (9) produces simpler equations.

For odd $\tau(p) = t = 4\,s \pm 1$, a reduction modulo 4 proves that the system of equations (10) has no solutions. Thus, suppose that $\tau(p) = t = 2\,s$ is even. Now it will be shown that this system of equations has finitely many integers solutions of the form $(p, X, Y) \in \mathbb{P}\times\mathbb{Z}\times\mathbb{Z}$. This is done in several steps.

Eliminating the prime variable $p$ yields an algebraic equation $F_t(X, Y) = 0$ over the integers, (use resultant or Grobner basis):

$$\begin{aligned} & 6\,973\,568\,802\,t^{11} - 3\,486\,784\,401\,t^{12} - 141\,021\,057\,996\,t^9\,X + 153\,418\,513\,644\,t^{10}\,X - \\ & 125\,352\,051\,552\,t^9\,X^2 + 136\,372\,012\,128\,t^{10}\,X^2 - 9\,415\,331\,872\,128\,t^7\,X^3 + \\ & 8\,535\,877\,796\,160\,t^8\,X^3 - 2\,690\,094\,820\,608\,t^7\,X^4 - 157\,154\,675\,198\,976\,t^5\,X^5 - \\ & 170\,970\,470\,820\,864\,t^6\,X^5 - 749\,464\,271\,142\,912\,t^3\,X^7 - 11\,414\,917\,360\,484\,352\,t^4\,X^7 - \\ & 148\,190\,367\,489\,196\,032\,t^2\,X^9 - 578\,220\,423\,796\,228\,096\,X^{11} - 10\,847\,773\,692\,t^9\,Y - \\ & 6\,198\,727\,824\,t^{10}\,Y + 1\,549\,681\,956\,t^{11}\,Y + 376\,056\,154\,656\,t^8\,X\,Y - \\ & 340\,930\,030\,320\,t^9\,X\,Y + 55\,712\,022\,912\,t^8\,X^2\,Y + 20\,922\,959\,715\,840\,t^6\,X^3\,Y - \\ & 15\,174\,893\,859\,840\,t^7\,X^3\,Y + 279\,386\,089\,242\,624\,t^4\,X^5\,Y + 227\,960\,627\,761\,152\,t^5\,X^5\,Y + \\ & 999\,285\,694\,857\,216\,t^2\,X^7\,Y + 10\,146\,593\,209\,319\,424\,t^3\,X^7\,Y + 65\,862\,385\,550\,753\,792\,t\,X^9\,Y + \\ & 28\,927\,396\,512\,t^8\,Y^2 + 1\,377\,495\,072\,t^9\,Y^2 - 417\,840\,171\,840\,t^7\,X\,Y^2 + \\ & 303\,048\,915\,840\,t^8\,X\,Y^2 - 18\,598\,186\,414\,080\,t^5\,X^3\,Y^2 + 10\,116\,595\,906\,560\,t^6\,X^3\,Y^2 - \\ & 186\,257\,392\,828\,416\,t^3\,X^5\,Y^2 - 101\,315\,834\,560\,512\,t^4\,X^5\,Y^2 - 444\,126\,975\,492\,096\,t\,X^7\,Y^2 - \\ & 2\,254\,798\,490\,959\,872\,t^2\,X^7\,Y^2 - 36\,427\,091\,904\,t^7\,Y^3 + 247\,608\,990\,720\,t^6\,X\,Y^3 - \\ & 134\,688\,407\,040\,t^7\,X\,Y^3 + 8\,265\,860\,628\,480\,t^4\,X^3\,Y^3 - 2\,997\,509\,898\,240\,t^5\,X^3\,Y^3 + \\ & 55\,187\,375\,652\,864\,t^2\,X^5\,Y^3 + 15\,009\,753\,268\,224\,t^3\,X^5\,Y^3 + 65\,796\,588\,961\,792\,X^7\,Y^3 + \\ & 28\,570\,268\,160\,t^6\,Y^4 - 82\,536\,330\,240\,t^5\,X\,Y^4 + 29\,930\,757\,120\,t^6\,X\,Y^4 - \\ & 1\,836\,857\,917\,440\,t^3\,X^3\,Y^4 + 333\,056\,655\,360\,t^4\,X^3\,Y^4 - 6\,131\,930\,628\,096\,t\,X^5\,Y^4 - \\ & 15\,237\,476\,352\,t^5\,Y^5 + 14\,673\,125\,376\,t^4\,X\,Y^5 - 2\,660\,511\,744\,t^5\,X\,Y^5 + \\ & 163\,276\,259\,328\,t^2\,X^3\,Y^5 + 5\,643\,509\,760\,t^4\,Y^6 - 1\,086\,898\,176\,t^3\,X\,Y^6 - \\ & 1\,433\,272\,320\,t^3\,Y^7 + 238\,878\,720\,t^2\,Y^8 - 23\,592\,960\,t\,Y^9 + 1\,048\,576\,Y^{10} = 0. \end{aligned} \qquad (11)$$

The system of partial derivative equations

$$\partial_X F_t(X, Y) = 0, \quad \partial_Y F_t(X, Y) = 0, \qquad (12)$$

does not have any integer roots $(X, Y)$, see [15, p. 150]. Moreover, it satisfies the congruences

$$F_t(X, Y) \not\equiv 0 \bmod X, \qquad F_t(X, Y) \not\equiv 0 \bmod Y. \qquad (13)$$

In particular, for any fixed nonzero even integer $t \ne 0$, the $x$-coordinate

$$4X = p^3(p^7 + 1) - pt, \tag{14}$$

and the $y$-coordinate cannot be written as a multiple

$$4Y = 13 p^4(p^7 + 1) - (22 p^2 - 9)t \ne f(p)X, \tag{15}$$

for some polynomial $f(x) \in \mathbb{Z}[x]$. For example,

$$\gcd(X, Y) = \gcd(13 p^4(p^7 + 1) - (22 p^2 - 9)t, \; p^3(p^7 + 1) - pt) = 1. \tag{16}$$

This means that $F_t(X, Y) \in \mathbb{Z}[X, Y]$ does not have a rational parametrization as a function of the prime variable $p \in \mathbb{P} = \{2, 3, 5, \ldots\}$. Thus, the algebraic equation $F_t(X, Y) = 0$ in (11) is nonsingular, see [15, p. 256].

By the Thue-Siegel theorem, it follows that it has finitely many integers solutions $(X, Y) \in \mathbb{Z} \times \mathbb{Z}$. Consequently, the subset of integer points is finite. That is,

$$\#\{(X, Y) \in \mathbb{Z} \times \mathbb{Z} : F_t(X, Y) = 0\} = O(1) \tag{17}$$

Ergo, for each fixed $0 \ne t \in \mathbb{Z}$, the cardinality of the subset of primes $\{p \in \mathbb{P} : \tau(p) = t \ne 0\}$ is finite:

$$\#\{p \in \mathbb{P} : \tau(p) = t \ne 0\} \ll \#\{(X, Y) \in \mathbb{Z} \times \mathbb{Z} : F_t(X, Y) = 0\}, \tag{18}$$

this follows from (17). ∎

Assuming similar congruences such as (9) can be determined for $\tau_k(n)$ the $n$th coefficients $\tau_k(n)$ of the discriminant functions $\Delta_k(z) = \sum_{n=1}^{\infty} \tau_k(n) q^n$ of level $N = 1$ and weights $k = 12, 16, 18, 20, 22,$ and $26$, the same analysis, mutatis mutandis, should be applicable.

These are special cases of the Lang-Trotter conjecture for cusp forms of level $N = 1$ and weights $k = 12, 16, 18, 20, 22,$ and $26$. Specifically, it claims that for any fixed integer $t \ne 0$, the subsets of primes

$$\{p \in \mathbb{P} : \tau_k(p) = t\} \tag{19}$$

is finite.

**Note 1.** More generally, the Thue-Siegel theorem states that a nonsingular algebraic curve $C : f(x, y) = 0$ such as (10) has finitely many integer solutions, that is,

$$\# C(\mathbb{Z}) = \#\{(x, y) \in \mathbb{Z} \times \mathbb{Z} : f(x, y) = 0\} = O(1). \tag{20}$$

The precise information on the torsion groups $C(\mathbb{Z})$ for algebraic curves of genus $g = 1$ are given in [27], et alii, some information on the ranks of such curves are given in [2], [26], and information on the sizes of the integral points are given in [29].

## 4. Proof Of Theorem 1

Let $\mathbb{P} = \{\ 2,\ 3,\ 5,\ ...\ \}$ be the set of primes. A subset of primes $\mathcal{P} = \{\ p \leqslant x : p \in \mathbb{P}\ \}$ generates a subset of integers

$$\{\ \tau(p) : p \leqslant x\ \} \tag{21}$$

of cardinality $\#\{\ \tau(p) : p \leqslant x\ \} \gg x / \log x$. This follows from Theorem 2, which proves that these subsets are finite. More generally, the multiplicative span of the cross product of $d \geqslant 1$ copies $\mathcal{P} \times \mathcal{P} \times \cdots \times \mathcal{P}$ generates a subset of integers

$$\{\ \tau(p_1\ p_2 \cdots p_d) : p_i \leqslant x\ \} \tag{22}$$

of cardinality $\#\{\ \tau(p_1\ p_2 \cdots p_d) : p_i \leqslant x\ \} \gg (x / \log x)^d$. For example, if $d \geqslant \log x$, then $\#\{\ \tau(p_1\ p_2 \cdots p_d) : p_i \leqslant x\ \} \gg x$.

## 5. The Subset $V_0$

The interesting case $t = 0$ was not considered in Section 2. In this case the system of equation (10) is singular, and it has infinitely many prime solutions $p \in \{\ p = 4m - 1 : m \in \mathbb{N}\ \} \subset \mathbb{P} = \{2,\ 3,\ 5,\ 7,\ ...\}$.

In fact, for $\tau(p) = t = 0$, the system of equations (10) reduces to the simpler form

$$\begin{aligned} p^3\left(p^7 + 1\right) - 4\,X &= 0, \\ 13\,p^4\left(p^7 + 1\right) - 4\,Y &= 0, \\ \tau(p) &= 0. \end{aligned} \tag{23}$$

The corresponding algebraic equation $F_0(X, Y) = 0$ is

$$-578\,220\,423\,796\,228\,096\,X^{11} + 65\,796\,588\,961\,792\,X^7\,Y^3 + 1\,048\,576\,Y^{10} = 0. \tag{24}$$

The system of partial derivative equations

$$\begin{aligned} \partial_X F_0(X, Y) &= -4\,902\,227\,890\,625\,X^6\left(15\,104\,375\,X^4 - 448\,Y^3\right) = 0, \\ \partial_Y F_0(X, Y) &= 320\,Y^2\left(2\,941\,336\,734\,375\,X^7 + 32\,768\,Y^7\right) = 0, \end{aligned} \tag{25}$$

has a double root at $(X, Y) = (0, 0)$, see [15, p. 150]. Moreover, it satisfies the congruences

$$F_0(X, Y) \equiv 0 \bmod X, \quad F_0(X, Y) \equiv 0 \bmod Y. \tag{26}$$

This means that it has a rational parametrization as a function of the primes $p \in \{ p = 4m - 1 : m \in \mathbb{N} \}$. In particular, the $x$-coordinate

$$4X = p^3(p^7 + 1), \tag{27}$$

and the $y$-coordinate can be written as a multiple

$$4Y = 13 p^4(p^7 + 1) = 13 p X. \tag{28}$$

Therefore, the algebraic equation $F_0(X, Y) = 0$ derived from the system of equations in (23) is singular, see [15, p. 256].

**Theorem 3.** Let $\tau(n)$ be the $n$th coefficient of the discriminat function $\Delta(z) = \sum_{n=1}^{\infty} \tau(n) q^n$. Then, $\tau(n) \neq 0$ for all integers $n \in \mathbb{N}$. In particular, the cardinality of the subset $V_0$ is finite: $\#V_0 = \#\{ p \in \mathbb{N} : \tau(p) = 0 \} = 0$.

**Proof:** The system of equations (23) has infinitely many prime solutions $p \in \{ p = 4m - 1 : m \in \mathbb{N} \}$. By Dirichlet Theorem, this is a subset of primes of positive density. But, this contradicts the well known zero density result

$$\#\{ \text{prime } p \in \mathbb{N} : \tau(p) = 0 \} = O\left(x \big/ (\log x)^{3/2 - \epsilon}\right), \tag{29}$$

see [14], [16], [21], [23], [24, Theorem 1.6], and similar references. The cardinalities of these subsets of primes are compactly compared in the inequalities

$$\frac{1}{2} \frac{x}{\log x} + o\left(\frac{x}{\log x}\right) = \#\{ \text{prime } p = 4m - 1 : m \in \mathbb{N} \}$$

$$\leq \#\{ \text{prime } p \in \mathbb{N} : \tau(p) = 0 \} \tag{30}$$

$$= O\left(\frac{x}{(\log x)^{3/2 - \epsilon}}\right),$$

where $\epsilon > 0$ is arbitrarily small. Clearly, this is a contradiction, and it implies that $\tau(p) \neq 0$ for all primes $p \geq 2$. ∎

## 6. The Subset $V_{pt}$

There is some interest in counting the subset of primes

$$\{ p \in \mathbb{P} : \tau(p) = pZ \neq 0, \; Z \in \mathbb{Z} \}. \tag{31}$$

The cardinality of this subset is unknown, it might be finite or infinite but has a very slow rate of growth. For a cusp form $f(z) = \sum_{n=1}^{\infty} \lambda(n) q^n$ of weight $k > 2$, it is expected that the cardinality of this subset of primes has a counting function such as

$$\pi(x, f) = \#\{ p \leq x : \lambda(p) = pZ \neq 0, \; Z \in \mathbb{Z} \} = O(\log\log x). \tag{32}$$

Confer [5] for discussions and numerical results for $\tau_k(p)$, $k = 12, 16, 18, 20, 22, 26$. These primes satisfy the congruence

$$\tau(p) \equiv 0 \bmod p. \tag{33}$$

Currently, the primes $p = 2, 3, 5, 7, 2411, 7758337633 \leq 10^{10}$ are known to satisfy this congruence, see [25] for other details.

***Theorem 4.*** Let $\tau(n)$ be the $n$th coefficient of the discriminant function $\Delta(z) = \sum_{n=1}^{\infty} \tau(n) q^n$. Then, the subset of integers $V_{pt} = \{ p \in \mathbb{N} : \tau(p) = pt, \; t \in \mathbb{Z} \}$ is finite.

***Proof:*** To derive some Diophantine equations associated with this case, let $\tau(p) = pZ$ for variable prime $p \in \mathbb{P}$, and variable integer $Z \in \mathbb{Z}$. Now consider

$$\begin{aligned} p^3(p^7+1) - p\,\tau(p) - 4X' &= 0, \\ 13 p^4(p^7+1) - (22 p^2 - 9)\tau(p) - 4Y' &= 0, \\ \tau(p) &= pZ, \end{aligned} \tag{34}$$

where $X', Y' \in \mathbb{Z}$ are integers variables, see (9). Any combination of the moduli specified by the system of congruences (9) can be selected; the choice of 4 for both congruences in (9) produces managable equations. To proceed, let $X = X'/p$, $Y = Y'/p$, and replace $\tau(p) = pZ$ in (34). Then, it reduces to

$$\begin{aligned} p^2(p^7+1) - pZ - 4X &= 0, \\ 13 p^3(p^7+1) - (22 p^2 - 9)Z - 4Y &= 0, \end{aligned} \tag{35}$$

where $\gcd(X, Y) = 1$. For odd $\tau(p) = pZ = p(4s \pm 1)$, a reduction modulo 4 proves that the systems of equations (34) and (35) have no solutions. Thus, suppose that $\tau(p) = pZ$ is even. Now it will be shown that this system of equations has finitely many integers solutions of the form $(p, X, Y, Z) \in \mathbb{P} \times \mathbb{Z} \times \mathbb{Z} \times \mathbb{Z}$. This is done in several steps.

Eliminating the prime variable $p$ yields an algebraic equation $F(X, Y, Z) = 0$ over the integers, (use resultant or

)

Grobner basis):

$$\begin{aligned}
&11\,119\,623\,534\,542\,848\,X^{10} - 16\,449\,147\,240\,448\,X^7\,Y^2 - 262\,144\,Y^9 + 74\,021\,162\,582\,016\,X^7\,Y\,Z - \\
&1\,118\,542\,012\,350\,464\,X^8\,Y\,Z + 1\,532\,982\,657\,024\,X^5\,Y^3\,Z + 5\,308\,416\,Y^8\,Z - \\
&83\,273\,807\,904\,768\,X^7\,Z^2 + 2\,516\,719\,527\,788\,544\,X^8\,Z^2 - 10\,347\,632\,934\,912\,X^5\,Y^2\,Z^2 + \\
&30\,440\,655\,618\,048\,X^6\,Y^2\,Z^2 - 40\,819\,064\,832\,X^3\,Y^4\,Z^2 - 47\,775\,744\,Y^7\,Z^2 + \\
&23\,282\,174\,103\,552\,X^5\,Y\,Z^3 - 136\,982\,950\,281\,216\,X^6\,Y\,Z^3 + 367\,371\,583\,488\,X^3\,Y^3\,Z^3 - \\
&29\,156\,474\,880\,X^4\,Y^3\,Z^3 + 271\,724\,544\,X\,Y^5\,Z^3 + 250\,822\,656\,Y^6\,Z^3 - \\
&17\,461\,630\,577\,664\,X^5\,Z^4 + 154\,105\,819\,066\,368\,X^6\,Z^4 - 1\,239\,879\,094\,272\,X^3\,Y^2\,Z^4 + \\
&196\,806\,205\,440\,X^4\,Y^2\,Z^4 - 3\,056\,901\,120\,X\,Y^4\,Z^4 - 4\,929\,859\,584\,X^2\,Y^4\,Z^4 - \\
&846\,526\,464\,Y^5\,Z^4 + 1\,859\,818\,641\,408\,X^3\,Y\,Z^5 - 442\,813\,962\,240\,X^4\,Y\,Z^5 + \\
&13\,756\,055\,040\,X\,Y^3\,Z^5 + 44\,368\,736\,256\,X^2\,Y^3\,Z^5 + 1\,904\,684\,544\,Y^4\,Z^5 + 13\,436\,928\,Y^5\,Z^5 - \\
&1\,046\,147\,985\,792\,X^3\,Z^6 + 332\,110\,471\,680\,X^4\,Z^6 - 30\,951\,123\,840\,X\,Y^2\,Z^6 - \\
&149\,744\,484\,864\,X^2\,Y^2\,Z^6 - 2\,857\,026\,816\,Y^3\,Z^6 - 151\,165\,440\,Y^4\,Z^6 + 51\,732\,592\,704\,X^3\,Z^7 + \\
&34\,820\,014\,320\,X\,Y\,Z^7 + 224\,616\,727\,296\,X^2\,Y\,Z^7 + 2\,678\,462\,640\,Y^2\,Z^7 + \\
&680\,244\,480\,Y^3\,Z^7 - 15\,669\,006\,444\,X\,Z^8 - 126\,346\,909\,104\,X^2\,Z^8 - 1\,205\,308\,188\,Y\,Z^8 - \\
&382\,637\,520\,X\,Y\,Z^8 - 1\,530\,550\,080\,Y^2\,Z^8 + 860\,934\,420\,X\,Z^9 - 15\,152\,445\,792\,X^2\,Z^9 + \\
&1\,721\,868\,840\,Y\,Z^9 - 774\,840\,978\,Z^{10} - 172\,186\,884\,Y\,Z^{10} + 387\,420\,489\,Z^{11} = 0\,.
\end{aligned} \quad (36)$$

The system of partial derivative equations

$$\partial_X F(X, Y, Z) = 0\,, \quad \partial_Y F(X, Y, Z) = 0\,, \quad \partial_Z F(X, Y, Z) = 0\,, \quad (37)$$

does not have any integer roots $(X, Y, Z)$. Thus, the algebraic equation $F(X, Y, Z) = 0$ in (36) is nonsingular, see [15, p. 150].

By the Thue-Siegel theorem, it follows that it has finitely many integers solutions $(X, Y, Z) \in \mathbb{Z} \times \mathbb{Z} \times \mathbb{Z}$. Consequently, the subset of integer points is finite. That is,

$$\#\{\,(X, Y, Z) \in \mathbb{Z} \times \mathbb{Z} \times \mathbb{Z} : F(X, Y, Z) = 0\,\} = O(1). \quad (38)$$

Ergo, the cardinality subset of the subset of primes $\{\,p \in \mathbb{P} : \tau(p) = pZ \neq 0,\ Z \in \mathbb{Z}\,\}$ is finite:

$$\#\{\,p \in \mathbb{P} : \tau(p) = pZ \neq 0,\ Z \in \mathbb{Z}\,\} \ll \#\{\,(X, Y, Z) \in \mathbb{Z} \times \mathbb{Z} \times \mathbb{Z} : F(X, Y, Z) = 0\,\}\,. \quad (39)$$

The inequality follows from (38). ∎